\documentclass[reqno,a4paper]{amsart}

\usepackage{amsfonts,amsmath,amssymb,amsthm}
\usepackage{alltt,esvect,graphicx,stmaryrd,tikz,xcolor}
\usepackage{pgffor}
\usepackage[arrow, matrix, curve]{xy}
\usepackage[titletoc]{appendix}
\usepackage[linesnumbered]{algorithm2e}
\usepackage[T1]{fontenc}
\usepackage[utf8]{inputenc}
\usepackage{helvet}
\usetikzlibrary{calc,intersections,through,shapes,arrows}
\usepackage{mathtools}
\usepackage{subfig}
\usepackage{url}
\usepackage{booktabs}
\usepackage{comment}
\usepackage[textsize=tiny]{todonotes}

\usepackage{enumitem} 
\usetikzlibrary{cd}
\usetikzlibrary{patterns}
\usetikzlibrary{math}
\usetikzlibrary{shapes}
\usetikzlibrary{decorations.fractals}

\usepackage{xstring}

\usepackage{cancel}
\usepackage{hyperref}

\newtheorem{theorem}{Theorem}[section]

\newtheorem{proposition}[theorem]{Proposition}
\newtheorem{corollary}[theorem]{Corollary}

\theoremstyle{definition}
\newtheorem{definition}[theorem]{Definition}
\newtheorem{example}[theorem]{Example}

\newtheorem{remark}[theorem]{Remark}

\definecolor{colorX}{RGB}{240,0,0}
\definecolor{colorY}{RGB}{0,0,160}
\definecolor{colorZ}{RGB}{0,160,0}

\newcommand{\Z}{\mathbb{Z}}

\newcommand{\R}{\mathbb{R}}
\newcommand{\C}{\mathbb{C}}

\newcommand{\T}{\mathbb{T}}

\renewcommand{\P}{\mathbb{P}}

\newcommand{\CL}{\mathcal{L}}

\newcommand{\CE}{\mathcal{E}}
\newcommand{\CH}{\mathcal{H}}
\newcommand{\CO}{\mathcal{O}}
\newcommand{\CP}{\mathcal{P}}
\newcommand{\CR}{\mathcal{R}}

\newcommand{\gH}{\operatorname{H}}

\newcommand{\rk}{\mathop{\rm rk}\nolimits}

\newcommand{\mb}{\mathbf}
\newcommand{\mc}{\mathcal}

\newcommand{\bp}{\begin{para}}
\newcommand{\ep}{\end{para}}
\newcommand{\bps}{\begin{paras}}
\newcommand{\eps}{\end{paras}}
\newcommand{\benum}{\begin{enumerate}[{\rm(i)}]}
\newcommand{\eenum}{\end{enumerate}}

\newcommand{\Div}{\operatorname{Div}}
\newcommand{\Cl}{\operatorname{Cl}}

\newcommand{\kss}{\scriptscriptstyle}
\newcommand{\kbb}{{\kss \bullet}}
 
\newcommand{\ko}{\overline}

\newcommand{\opn}{\operatorname}

\newcommand{\Nef}{\operatorname{Nef}}

\newcommand{\TX}{\mb{T}}

\newcommand{\ssect}[1]{Subsection~\ref{#1}}

\newcommand{\CD}{{\mathcal D}}
\newcommand{\CT}{{\mathcal T}}

\newcommand{\drawExSet}[1]{
 \StrSubstitute{#1}{,}{/}[\replacedSlash]
 \StrSubstitute{\replacedSlash}{;}{,}[\xycoords]
 \foreach [count=\i] \x/\y in \xycoords
 {
  \ifthenelse {1=\i}
   {\fill (\x,\y) circle (5pt); \draw (\x,\y) circle (7pt);}
   {\fill (\x,\y) circle (5pt);}
 }
}
\definecolor{oliwkowy}{HTML}{627037}
\definecolor{lightblue}{RGB}{135,206,250}
\definecolor{darkblue}{RGB}{0,0,160}
\definecolor{darkgreen}{RGB}{0,160,0}
\definecolor{veryPeri}{RGB}{102,103,171}
\definecolor{intOrange}{rgb}{1.0,.310,.0}
\definecolor{MidnightBlue}{RGB}{102,103,171}

\definecolor{cklaus}{rgb}{0.0,.288,.378}
\definecolor{candreas}{RGB}{10,100,20}

\newcounter{para}[section]
\newenvironment*{para}[1]{\refstepcounter{para}\noindent\ignorespaces{\bf\thepara.~#1.}}{\ignorespacesafterend\bigskip}
\newenvironment*{paras}[1]{{\bf #1.}}{\ignorespacesafterend\bigskip}
\numberwithin{para}{section}

\setlist[enumerate,1]{label = (\roman*),ref = \theenumii.\roman*}

\newcommand{\orb}{\operatorname{orb}}

\newcommand{\into}{\hookrightarrow}

\newif\ifnotext

\definecolor{colorAB}{RGB}{240,0,0}
\definecolor{colorBC}{RGB}{0,0,160}
\definecolor{colorCD}{rgb}{1.0,.310,.0}
\definecolor{colorDA}{RGB}{0,160,0}

\newcommand{\TT}{{\opn{T}}}

\newcommand{\Strat}{{\mathcal S}}

\newcommand{\Pol}{\operatorname{Pol}}

\definecolor{cRayA}{rgb}{0.0,0.6,0.0}
\definecolor{cRayB}{rgb}{0.0,0.0,1.0}
\definecolor{cRayC}{rgb}{1.0,0.0,0.0}

\definecolor{cConeA}{rgb}{0.0,0.6,0.0}
\definecolor{cConeB}{rgb}{0.0,0.0,1.0}
\definecolor{cConeC}{rgb}{1.0,0.0,0.0}

\newif\iftext

\renewcommand{\S}{\mathbb{S}}
\definecolor{green}{rgb}{0.0,0.6,0.0}
\definecolor{blue}{rgb}{0.0,0.0,1.0}
\definecolor{red}{rgb}{1.0,0.0,0.0}

\renewcommand{\Z}{\mathbb{Z}}
\newcommand{\NE}{\operatorname{NE}}
\newcommand{\pc}{{\CP}}
\newcommand{\foc}{{\sigma}} 
\newcommand{\ffoc}{{f}}
\newcommand{\pRel}{{\CR}}
\newcommand{\face}{\operatorname{face}}

\begin{document}
\parindent0mm

\title[Acyclic toric sheaves]
{Acyclic toric sheaves}

\author[K.~Altmann]{Klaus Altmann}
\address{Institut f\"ur Mathematik,
FU Berlin,
K\"onigin-Luise-Str.~24-26,
D-14195 Berlin
}
\email{altmann@math.fu-berlin.de}
\author[A.~Hochenegger]{Andreas Hochenegger}
\address{
Dipartimento di Matematica ``Francesco Brioschi'',
Politecnico di Milano,
via Bonardi 9,
20133 Milano 
}
\email{andreas.hochenegger@polimi.it}
\author[F.~Witt]{Frederik Witt}
\address{
Fachbereich Mathematik,
U Stuttgart,
Pfaffenwaldring 57,
D-70569 Stuttgart
}
\email{witt@mathematik.uni-stuttgart.de}

\thanks{  
\hfill\newline
}

\begin{abstract}
Let $\CE$ be a torus-linearised reflexive sheaf over a smooth projective toric variety. Generalising a theorem of Perlman and Smith,
we prove an explicit sufficient condition for $\CE$ to be acyclic via 
Weil decorations.
\end{abstract}

\maketitle

\section{Introduction}
\label{sec:Intro}
For a given sheaf it is natural to ask whether it is {\em acyclic}, that is, it has no higher cohomology except possibly for $\opn{H}^0$, or even {\em immaculate}, that is, it has no cohomology at all. In the setting of toric geometry, Perlman and Smith~\cite{MFO} give a neat condition for the acyclity of {\em toric vector bundles}, that is, locally free sheaves with a linearised torus action.

\subsection{The Perlman-Smith theorem.}
\label{subsec:PST}
To state it, we let $X$ be a $d$-dimensional, smooth and projective toric variety over $\C$ with embedded torus $j_\TT\colon\TT\into X$ and fan $\Sigma$. As usual, $M$ denotes the character lattice of $\TT$ and $N$ its dual. 
Following Klyachko~\cite{klyachko}, 
a toric vector bundle $\CE$
can be described in terms of $\Z$-descending filtrations of the $\C$-vector space $E$ given by the torus-invariant sections in $\Gamma(\TT,\CE)$, that is,
\[
E_\rho^\kbb\;=\;[\ldots\supseteq E_\rho^{\ell-1}\supseteq E_\rho^{\ell}\supseteq E_\rho^{\ell+1}\supseteq\ldots]
\]
with $E_\rho^\ell=E$ for $\ell\ll0$ and $=0$ for $\ell\gg0$. They are parametrised by the one-dimensional cones or {\em rays} $\rho\in\Sigma(1)$ of the fan. This datum is supposed to satisfy a further compatibility condition, see for instance also~\cite{payne} for a short and concise statement. The filtrations induce the quantities 
\[
\lambda_\rho:=\max\{\ell\in\Z\mid E_\rho^\ell\neq 0\}\quad\mbox{and}\quad\mu_\rho:=\max\{\ell\in\Z\mid E_\rho^\ell=E\}.
\]

\medskip

Finally, we recall that a {\em primitive collection} $\pc$ is a subset of 
rays $\{\rho_1,\ldots,\rho_k\}$ defining a minimal non-face, that is, 
every strict subset of $\pc$ is a set of rays $\tau(1)$ for some cone 
$\tau\in\Sigma$, but $\pc$ itself is not. We let $\foc(\pc)$ be the uniquely determined cone in $\Sigma$ for which $\sum_{\rho\in\pc}\rho$ 
lies in the interior; we note in passing that we always identify a ray with its primitive generator. Consider the natural map $\pi\colon\Z^{\Sigma(1)}\to N$ that assigns the base vector $e_\rho$ of $\Z^{\Sigma(1)}$ to $\rho\in N$. The {\em focus of $\pc$} is the unique element
\[
\ffoc(\pc):=\sum_{\rho\in\foc(\pc)(1)}\ffoc_\rho\cdot e_\rho
\in\Z_{\geq 1}^{\foc(\pc)(1)} 
\]
such that the so-called {\em primitive relation}
\[
\pRel(\pc):=\sum_{\rho\in\pc}e_\rho-\ffoc(\pc)
\]
lies in the kernel of $\pi$ which is $\Cl(X)^\vee$, the dual of the class group of $X$.

\medskip

{\bf Theorem (Perlman-Smith~\cite{MFO}).}
{\em Let $\CE$ be a toric vector bundle on $X$ and $D=\sum_{\rho\in\Sigma(1)}a_\rho D_\rho$ be a $\TT$-invariant divisor such that for all primitive collections $\pc\subseteq\Sigma(1)$, the inequality
\begin{equation}
\label{eq:PS}
\sum_{\rho\in \pc} (a_{\rho}+\mu_{\rho}) \geq
\sum_{\rho\in \foc(\pc)(1)} \ffoc_\rho\cdot(a_{\rho}+\lambda_{\rho})
\end{equation}
holds for the focus $\ffoc(\pc)=\sum_{\rho\in\foc(\pc)(1)}\ffoc_\rho\cdot\rho$ of $\pc$. Then $\CE(D)=\CE\otimes\CO_X(D)$ is acyclic, that is, $\opn{H}^i(X,\CE(D))=0$ for all $i\geq 1$.
}

\begin{remark}
The notion of a primitive collection goes back to Batyrev~\cite{picRank3}. 
As he showed and we use below, the set of primitive relations 
$\pRel(\pc)$ coming from the primitive collections spans the Mori cone 
$\NE(X)\subseteq 
\Cl(X)^\vee \otimes_\Z \R$
which is dual to the nef cone $\Nef(X)\subseteq
\Cl(X)\otimes_\Z\R$.
\end{remark}

\subsection{Weil decorations}
So-called {\em Weil decorations} were introduced in the preprint \cite{TRS} and provide an alternative description of toric vector bundles and in fact of toric sheaves, that is, reflexive sheaves with a $\TT$-linearised action. Roughly speaking, a Weil decoration of $\CE$ consists of a finite collection of $\TT$-invariant divisors $\{D_1,\ldots,D_n\}$ relating to the vector space $E$ from \ssect{subsec:PST}. From~\cite{TRS}, we immediately derive the following vanishing theorem and its corollary, see~\ssect{subsec:CanRes}.

\begin{theorem}
\label{thm:VanThm}
If $k_0$ is an integer such that $\gH^k\!\big(X,\CO(D_i)\big)=0$ for $i=1,\ldots,n$ and all $k\geq k_0$, then also $\gH^k(X,\CE)=0$ for all $k\geq k_0$. In particular, this implies that if 
\begin{enumerate}
\item[{\rm (i)}] all $\CO_X(D_i)$ are acyclic, then so is $\CE$. 

\smallskip

\item[{\rm (ii)}] all $\CO_X(D_i)$ are immaculate, then so is $\CE$. 
\end{enumerate}
\end{theorem}

\begin{remark}
\label{rem:LBPol}
Acyclicity or immaculacy of a line bundle $\CO_X(D)$ can be easily read off the divisor by writing $D$ as a (non-unique) difference $D_+-D_-$ of two torus-invariant nef Cartier divisors. These correspond to lattice polytopes $\nabla_+$ and $\nabla_-$ sitting inside $M_\R$. By~\cite{immaculate},~\cite{dop} we find
\[
\gH^k\!\big(X,\CO_X(D)\big)_m=\widetilde{\opn H}^{\raisebox{-3pt}{\scriptsize$k\!-\!1$}}\big(\nabla_-\setminus(\nabla_+-m)\big),
\]
for the cohomology of $\CO_X(D)$ in degree $m\in M$, where $\widetilde{\opn H}^{\raisebox{-3pt}{\scriptsize$k\!-\!1$}}$ denotes the $(k-1)$-th reduced singular cohomology with complex coefficients. In particular, we recover 
{\em Demazure vanishing}. Namely, a nef line bundle $\CO_X(D)$ is acyclic as follows from taking $D_-=0$.  
\end{remark}

\begin{corollary}
\label{coro:Neflic0}
If the divisors $D_1,\ldots,D_n$ are nef, then $\CE$ is acyclic. 
\end{corollary}

\begin{remark}
Since ample toric vector bundles are not necessarily acyclic, see for instance~\cite[Example 4.10]{HMP10} or~\cite[Remark 6.4]{TRS}, the corollary can be regarded as a substitute for Demazure vanishing of a line bundle, at least if $X$ is smooth and projective.
\end{remark}

Corollary~\ref{coro:Neflic0} also yields a straightforward proof for the following generalisation of the Perlman-Smith theorem, see \ssect{subsect:Proof}.

\begin{theorem}
\label{thm:Neflic}
Let $\CE$ be a toric sheaf and $D=\sum_{\rho\in\Sigma(1)}a_\rho D_\rho$ be a $\TT$-invariant divisor such that Inequality~\eqref{eq:PS} holds for all
primitive collections defining an extremal ray of the Mori cone. Then $\CE(D)$ is acyclic.
\end{theorem}

We shall give a geometric interpretation of the assumption of Theorem~\ref{thm:Neflic} in~\ssect{subsect:GeoInt}.

\begin{remark}
Theorem~\ref{thm:Neflic} will follow from the first implication in
\begin{center}
\eqref{eq:PS} $\Longrightarrow$ nefly decorated $\Longrightarrow$ acyclicly decorated $\Longrightarrow$ acyclic.
\end{center}
The converses are all wrong in general. The twisted tangent bundle $\CT_{\P^2}(-4)$ is acyclic with non-acyclic Weil decoration~\cite[Section 8.5]{TRS}; $\CO_{\P^1}\oplus\CO_{\P^1}(D_1-D_2)$ is acyclicly decorated, but not nefly~\cite[Example 3.9]{TRS}. In Example~\ref{exam:NefNotPS} we will construct an example of a nefly decorated toric vector bundle which doesn't satisfy~\eqref{eq:PS}.
\end{remark}

\subsection*{Acknowledgements}
The authors would like to thank Nathan Ilten for very helpful conversations.

\section{Weil decorations}
\label{sec:WeilDeco}
We first recall some technical background from~\cite{TRS}. In addition to the notation introduced in Section~\ref{sec:Intro}, we let $\Div_\TT(X)$ be the group of {\em torus-invariant Weil divisors}, that is, the free abelian group generated by the orbit closures
\[
D_\rho=\overline{\opn{orb}(\rho)},\quad\rho\in\Sigma(1).
\]
For us, a divisor will always mean a torus-invariant Weil divisor, i.e., a divisor in $\Div_\TT(X)$. Since $X$ is smooth, divisors $D$ stand in 1-1 correspondence with invertible subsheaves $\CL=\CO_X(D)$ of $j_{\TT*}\CO_\TT=\C[M]$ (we slightly abuse notation and drop here and in the sequel the sheafification of $\C[M]$).

\subsection{Polytopes and Weil divisors}
\label{subsec:PolWD}
A lattice polytope $\nabla$ is {\em compatible with} $\Sigma$ if its normal fan $\mc N(\nabla)$ is refined by $\Sigma$. With Minkowski sum as addition,
\[
\Pol^+(\Sigma):=\text{the set of compatible lattice polytopes}
\]
becomes a cancellative 
monoid
whose associated Grothendieck group is
\[
\Pol(\Sigma):=\text{the group of {\em virtual} polytopes}.  
\]
We think of the latter as formal differences $\nabla=\nabla_+-\nabla_-$ with $\nabla_\pm\in\Pol^+(\Sigma)$.

\medskip

Compatibility entails that for every cone $\sigma\in\Sigma$ there is a unique face of $\nabla\in\Pol^+(\Sigma)$ determined by
\begin{equation}
\label{eq:NablaFace}
\langle\face(\nabla,\sigma),a\rangle=\min\langle\nabla,a\rangle\quad\text{for all }a\in\sigma. 
\end{equation}
If $\sigma$ is full-dimensional, then $\face\big(\nabla,\sigma\big)$ is a vertex that we denote $\nabla(\sigma)$. The induced monomials $x^{\nabla(\sigma)}\in\C[M]$ define the invertible subsheaf $\CO_X(\nabla)$ of $(j_\TT)_*\CO_\TT$ by
\[
\CO_X(\nabla)|_{U_\sigma}:=x^{\nabla(\sigma)}\cdot\CO_{U_\sigma}=\C[\nabla_\sigma\cap M]\subseteq\C[M]
\]
where we used the {\em local polyhedron} $\nabla_\sigma:=\nabla(\sigma)+\sigma^\vee=\nabla+\sigma^\vee$. Furthermore,
\[
x^{\nabla(\sigma)}\in\gH^0(X,\CO_X(\nabla))=\C[\nabla\cap M],
\]
which implies that $\CO_X(\nabla)$ is globally generated. Equivalently, its associated divisor
\[
D_\nabla=-\sum_{\rho\in\Sigma(1)}\min\langle\nabla,\rho\rangle\cdot D_\rho
\]
is basepoint free, which for toric varieties is equivalent to being nef. 

\medskip

Conversely, starting with a divisor $D$, we consider the (possibly empty) {\em section polyhedron} 
\[
\gH^0_\R\!\big(\hspace{-5pt}\sum_{\rho\in\Sigma(1)}\!\!a_\rho D_\rho\big):=\big\{u\in M_\R\mid\langle u,\rho\rangle\geq-a_\rho,\,\rho\in\Sigma(1)\big\};   
\]
the integral points $\gH^0_\R(D)\cap M$ correspond to a $\C$-basis 
$\{x^m\}$ of global sections of $\CO_X(D)$. Then, $\gH^0_\R(D_\nabla)=\nabla$ for $\nabla\in\Pol^+(\Sigma)$. In particular, the restriction of the map $\gH^0_\R$ to nef divisors induces a 
monoid
isomorphism onto $\Pol^+(\Sigma)$. By projectivity, the Grothendieck group of $\opn{Nef}(\Sigma)$ is $\Div_\TT(X)$ so that $\gH^0_\R|_{\opn{Nef}}$ extends to a group isomorphism 
\[
\opn{P}\colon\Div_\TT(X)\to\Pol(\Sigma).
\]

\subsection{The lattice structure}
\label{subsec:lattice-structure}
The usual poset relation 
\[
D\leq D'\quad\text{if and only if}\quad D'-D\text{ is effective}
\]
turns $\Div_\TT(X)$ into a {\em lattice} with {\em meet}
\[
\min\!\big(\sum a_\rho D_\rho,\sum a'_\rho D_\rho\big):=\sum\min(a_\rho,a'_\rho)D_\rho.
\]
and {\em join} 
\[
\max\!\big(\sum a_\rho D_\rho,\sum a'_\rho D_\rho\big):=\sum\max(a_\rho,a'_\rho)D_\rho.
\]
By~\cite[Subsection 2.3.2]{TRS}, we have
\[
D\leq D'\quad\text{if and only if}\quad\opn{P}(D)\subseteq\opn{P}(D')
\]
for {\em virtual inclusion} $\subseteq$; it is honest inclusion for polytopes in $\Pol^+(\Sigma)$. Furthermore,
\[
\CO_X(D\wedge D')=\CO_X(D)\cap\CO_X(D')\subseteq\C[M],
\]
so we define the {\em virtual intersection} in $\Pol(\Sigma)$ by
\[
\nabla\Cap\nabla':=\opn{P}\!\big(\!\min(D_\nabla,D_{\nabla'})\big).
\]
For the join in $\Div_\TT(X)$, we have 
\[
\CO_X(D\vee D')=\big(\CO_X(D)+\CO_X(D')\big)^{\vee\vee},
\]
for $\CO_X(D\vee D')$ is the smallest reflexive subsheaf of $K(X)$ containing both $\CO_X(D)$ and $\CO_X(D')$, which by definition is the double dual of $\CO_X(D)+\CO_X(D')$. We therefore define the {\em virtual union} of two virtual polyhedra in $\Pol(\Sigma)$ by
\[
\nabla\Cup\nabla':=\opn{P}\!\big(\!\max(D_\nabla,D_{\nabla'})\big).
\]
In particular, we obtain the

\begin{proposition}
$\opn{P}\colon(\Div_\TT(X),\leq,\min,\max)\to(\Pol(\Sigma),\subseteq,\Cap,\Cup)$ defines a lattice isomorphism.
\end{proposition}

Table~\ref{table:DivPol} summarises the relationship between divisors, invertible sheaves and virtual polyhedra.
\begin{table}[ht]
\begin{tabular}{c|c|c}
$\opn{Div}(\Sigma)$ & invertible sheaves $\subseteq(j_\TX)_*\CO_\TX=\C[M]$ & $\Pol(\Sigma)$\\\hline\vspace{-5pt}&&\\
$D+D'$&$\CO_X(D)\cdot\CO_X(D')=\CO_X(D)\otimes\CO_X(D')$&$\opn{P}(D)+\opn{P}(D')$\\[3pt]
$-D$&$\CO_X(D)^{-1}=\mc{H}om_{\CO_X}\!\big(\CO_X(D),\CO_X\big)$&$-\opn{P}(D)$\\[3pt]
$D\leq D'$&$\CO_X(D)\subseteq\CO_X(D')$&$\opn{P}(D)\subseteq\opn{P}(D')$\\[3pt]
$\min(D,D')$&$\CO_X(D)\cap\CO_X(D')$& $\opn{P}(D)\Cap\opn{P}(D')$\\[2pt]
$\max(D,D')$&$\big(\CO_X(D)+\CO_X(D')\big)^{\vee\vee}$& $\opn{P}(D)\Cup\opn{P}(D')$\\
\end{tabular}
\caption{}
\label{table:DivPol}
\end{table}

\begin{remark}
Unlike $\subseteq$, $\Cap$ and $\Cup$ do {\em not} restrict to the usual union and intersection $\cup$ and $\cap$ for nef polytopes. 

\medskip

For instance, take the two ample polytopes $\nabla$ and $\nabla'$ for the del Pezzo surface obtained by blowing up $\P^2$ in two points fixed by the torus, see Figure~\ref{fig:PolOps} below. They intersect in the shaded square determined by the origin and the point $[2,2]$. The $D_{\rho_4}$-coefficient of its corresponding divisor is $4$ while the $D_{\rho_4}$-coefficient of the divisors $D_\nabla$ and $D_{\nabla'}$ is $5$. Hence $D_{\nabla\cap\nabla'}\not=\min(D_\nabla,D_{\nabla'})$. Note, however, that making $\nabla$ and $\nabla'$ sufficiently ample, $\nabla\Cap\nabla'$ actually becomes the honest intersection $\nabla\cap\nabla'$.

\medskip

Similarly, we see that the actual union of $\nabla$ and $\nabla'$ misses out the corners marked in orange, so again $\nabla\Cup\nabla'\not=\nabla\cup\nabla'$. Even worse, this inequality persists under making $\nabla$ and $\nabla'$ arbitrarily ample as the union is not even a polytope.

\begin{figure}[ht]
\begin{tikzpicture}[scale=0.3]
\draw[thick, color=black]
  (-4,0) -- (4,0) (0,-4) -- (0,4) (0,0) -- (-3.0,-3.0);
\draw[thick, color=black]
  (5,0) node{$\rho_1$} (0,5) node{$\rho_2$} (-5.5,0) node[left]{$\rho_3$} (-4.2,-3.8) node{$\rho_4$} (0,-5) node{$\rho_5$};
\end{tikzpicture}
\hspace{20pt}
\begin{tikzpicture}[scale=0.55]
\draw[color=oliwkowy!40] (-1.3,-0.3) grid (7.3,8.3);
\fill[pattern color=orange, pattern=north west lines]
  (2,3) -- (4,3) -- (4,5) -- (2,5) -- cycle;
\draw[very thick, color=green]
  (2,1) -- (2,7) -- (3,7) -- (4,6) --(4,1) -- cycle;
\draw[thick, color=green]
  (2.3,0.3) node{$\nabla'$};
\draw[very thick, color=blue]
  (0,3) -- (6,3) -- (6,4) -- (5,5) -- (0,5) -- cycle;
\draw[very thick, color=blue]
  (-0.4,2.3) node{$\nabla$};
\draw[very thick, color=magenta, dotted]
  (7.3,2.7) -- (1.7,8.3); 
\draw[very thick, color=magenta, dotted]
  (7.3,1.7) -- (0.7,8.3); 
\fill[color=red]
  (2,3) circle (6pt);
\fill[color=orange]
  (0,1) circle (4pt) (0,7) circle (4pt) (6,1) circle (4pt);
\end{tikzpicture}
\caption{Left hand side: The fan of the del Pezzo surface. Right-hand side: The red dot indicates the origin in $M_\R$ and fixes the position of the polytopes.}
\label{fig:PolOps}
\end{figure}
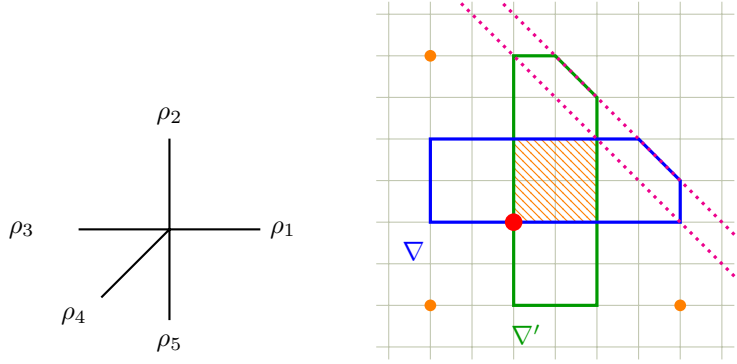
\end{remark}

\subsection{Weil decorations}
\label{subsec:WeilDeco}
We recall Weil decorations and some of their properties from~\cite{TRS}.

\begin{definition}
\label{def:WeilDec}
A {\em Weil decoration} on a $\C$-vector space $E$ is a map
\[
\CD\colon\fbox{$E^\times:=E\setminus\{0\}$}\to\Div(\Sigma) 
\]
which factorises over the projectivisation $\P(E)$ of $E$ and such that for all $e$, $e'\in E^\times$ with $e+e'\not=0$, the inequality
\begin{equation}
\label{eq:Basic} 
\CD(e+e')\geq\min\!\big(\CD(e),\CD(e')\big)
\end{equation}
holds true.
\end{definition}

One can show that~\eqref{eq:Basic} implies that the image of $\CD$ is finite 
(see \cite[Proposition 3.2]{TRS} for a proof). Therefore, it
has a lower bound which we denote $\CD(\eta)$. 
The notation hints at the fact that $\CD(e)=\CD(\eta)$ for generic $e\in E^\times$. Furthermore, the divisor
\[
\widehat D:=\bigvee\{\CD(e)\mid e\in E^\times\}.
\]
is well-defined.

\begin{remark}
In~\cite{TRS} we extended the Weil decoration to all of $E$ by setting formally $\CD(0)=\sum\infty\cdot D_\rho$ and working with the usual arithmetic in $\Z\cup\{\infty\}$. While this is convenient for the course of many proofs, it is often more natural to work over $E^\times$. At any rate, we can replace $(\infty)$ by any upper bound in $\Div_T(X)$, for instance by $\widehat D$ in view of the Perlman-Smith theorem.  
\end{remark}

Given a toric sheaf $\CE$ on $X$ we can build its {\em associated Weil decoration} as follows. By torsionfreeness of reflexive sheaves, the natural morphism $\CE\to(j_\TX)_*j_\TX^*\CE=(j_\TX)_*\CE|_\TX$ is actually injective. Let 
\[
E:=\text{ the torus invariant subspace of }\CE(\TX)=\Gamma(\TX,\CE);
\]
in particular, $\CE|_\TX$ is the sheafification of $E\otimes\C[M]$. As for rank $1$ sheaves we therefore consider 
\[
\CE\subseteq(j_\TX)_*\CE|_\TX=E\otimes(j_\TX)_*\CO_\TX=E\otimes\C[M]
\]
as a subsheaf of $E\otimes\C[M]$. For $e\in E^\times$ we get the saturated, hence reflexive subsheaf
\begin{equation}
\label{eq:CEe}
\CE(e):=\CE\cap\big(\C\cdot e \otimes(j_\TX)_*\CO_\TX\big)
\end{equation}
which means that $\CE(e)(U_\sigma)=\CE(U_\sigma)\cap(\C\cdot e\otimes\C[M])$. The isomorphic subsheaf $\CL(e)$ of $(j_\TX)_*\CO_\TX=\C[M]$ resulting via
\[
\begin{tikzcd}
\CE(e)\ar[r,hookrightarrow]\ar[d,"\cdot 1/e"',"\cong"]& e\cdot\C[M]\ar[d,"\cdot 1/e","\cong"']\\
\CL(e)\ar[r,hookrightarrow]&\C[M]
\end{tikzcd}
\]
induces a well-defined Weil divisor $D(e)$ with $\CO_X(D(e))=\CL(e)$. Furthermore, the assignment
\[
\CD_\CE\colon E^\times\to\Div(X),\quad\CD_\CE(e):=D(e)
\]
defines a Weil decoration~\cite[Subsection 3.2]{TRS}.

\begin{remark}
The Weil decoration $\CD_\CE$ actually determines $\CE$. For instance, we recover the Klyachko filtrations via 
\[
E^\ell_\rho=\{e\in E^\times\mid\,\big(\CD_\CE(e)\big)_\rho:=\text{the coefficient of }D_\rho\text{ in }\CD_\CE(e)\text{ is }\geq\ell\}\cup\{0\}.
\]
\end{remark}

\begin{remark}
\label{rem:DetaDNull}
In Section~\ref{sec:Intro} we defined the quantities
\[
\lambda_\rho=\max\{\ell\in\Z\mid E_\rho^\ell\neq 0\}\geq\mu_\rho:=\max\{\ell\in\Z\mid E_\rho^\ell=E\}
\]
which are determined by the Klyachko filtrations. Then 
\[
\CD(\eta)=\min\{\CD(e)\mid e\in E^\times\}=\!\sum_{\rho\in\Sigma(1)}\min\{\CD(e)_\rho\mid e\in E^\times\}\cdot D_\rho=\!\sum_{\rho\in\Sigma(1)}\mu_\rho\cdot D_\rho 
\]
and 
\[
\widehat D=\bigvee\big\{\CD(e)\mid e\in E^\times\big\}=\!\sum_{\rho\in\Sigma(1)}\lambda_\rho\cdot D_\rho.
\]
\end{remark}

\medskip

Since we are free to identify $\Div_\TT(X)$ with $\Pol(\Sigma)$ we will tacitly think of Weil decorations as being $\Pol(\Sigma)$-valued.

\begin{definition}
A toric sheaf $\CE$ is {\em nefly decorated} if the image of $\CD_\CE$ lies in $\Pol^+(\Sigma)$. 
\end{definition}

To keep notation tight we will usually drop the index $\CE$ and simply write $\CD$ instead of $\CD_\CE$. Pooling together the vectors $e\in E^\times$ into $\S(D)=\{e\in E^\times\mid\CD(e)=D\}$ and setting $\Strat:=\{\S(D)\mid D\in\CD(E^\times)\}$ yields a stratification $\Strat\cup\{0\}$ of $E$ by linear subspaces $\overline\S$ for $\S\in\Strat$. Further, $\Strat$ carries the natural partial order
\begin{equation}
\label{eq:StratPOS}
\S\leq\S'\quad\text{if and only if}\quad\S\subseteq\overline{\S'}.
\end{equation} 
The generic stratum whose divisor is the lower bound of $\CD(E)$ is denoted by $\eta$. The minimal stratum whose closure contains $\S(D)$ and $\S(D')$ defines a join $\S(D)\vee\S(D')$ on $(\Strat,\leq)$ and one easily verifies that
\[
\S(D)\vee\S(D')=\S(\min(D,D')). 
\]
Consequently, the Weil decoration $\CD$ descends to an anti-semilattice isomorphism onto the image of 
\[
\CD\colon(\Strat,\leq,\vee)\to(\Div(\Sigma),\leq,\min),\quad\S(D)\mapsto D, 
\]
that is,
\[
\CD(\S\vee\S')=\min\big(\CD(\S),\CD(\S')\big)
\]
for all strata $\S$, $\S'\in\Strat$.

\subsection{The canonical resolution}
\label{subsec:CanRes}
By \cite{TRS}, there is a particular resolution of toric sheaves by invertible ones which turns out to be useful here.
To define our complex we consider for $\S$, $\T\in\Strat$ the totally split sheaf
\[
\CE_{\S,\T}:=\overline\S\otimes_\C\CO_X(\T),
\]
where $\overline\S$ denotes closure of $\S$ in $E$ and $\CO_X(\T)$ the line bundle given by the divisor $\CD(\T)$. For $\CE_{\S,\S}$ we simply write $\CE_\S$. Further, we let
\[
\opn{ch}_\ell(\S,\T):=
\{\S=\S_0<\S_1<\ldots<\S_{\ell-1}<\S_\ell=\T\mid \S_1,\ldots,\S_{\ell-1}
\in\Strat\}
\]
be the set of strict chains of length $\ell$ in $\Strat$ starting at $\S_0=\S$ and terminating at $\S_\ell=\T$. Note that $\opn{ch}_\ell(\S,\T)\not=\varnothing$ implies $\S\leq\T$. For instance, 
\[
\opn{ch}_0(\S,\T)=\begin{cases}
\{\S\}&\text{if }\S=\T\\
\varnothing&\text{if }\S<\T
\end{cases}
\quad\text{and}\quad
\opn{ch}_1(\S,\T)=\begin{cases}
\varnothing&\text{if }\S=\T\\
\{\S<\T\}&\text{if }\S<\T
\end{cases}
. 
\]

\begin{proposition}[Corollary 6.10 in \cite{TRS}]
\label{prop:CanSeq}
For any toric sheaf $\CE$ of rank $r$, the complex
\begin{equation}
\label{eq:ResLB}
\fbox{$0\to\bigoplus_{\S\leq\T}\bigoplus_{\opn{ch}_r(\S,\T)}\!\CE_{\S,\T}\to\cdots\to\bigoplus_{\S\leq\T}\bigoplus_{\opn{ch}_0(\S,\T)}\!\CE_{\S,\T}$}\to\CE\to0
\end{equation}
with the usual differential is exact. In particular, the line bundles occurring in this complex are nef if $\CE$ is nefly decorated.
\end{proposition}

For the summand $\CE_{\S_0,\S_\ell}\subseteq\bigoplus_{\S\leq\T}\bigoplus_{\opn{ch}_\ell(\S,\T)}\CE_{\S,\T}$ corresponding to the chain $\S_0<\ldots<\S_\ell$, we could have written $\CE_{\S_0,\ldots,\S_\ell}$ for $\CE_{\S_0,\S_\ell}$ even though the summands do not depend on the inner terms. 

\medskip

The boxed complex in~\eqref{eq:ResLB} is thus quasi-isomorphic to $\CE$. Taking the standard spectral sequence of hypercohomology, namely,
\[
E^{-\ell,q}_1=\gH^q\big(X,\bigoplus_{\S\leq\T}\bigoplus_{\opn{ch}_\ell(\S,\T)}\!\CE_{\S,\T}\big)\Rightarrow\gH^{q-\ell}(X,\CE)
\]
for $0\leq\ell\leq r=\rk(\CE)$ and $0\leq q\leq d=\dim(X)$, yields the

\begin{theorem}
\label{thm:VanThm2}
If $k_0$ is an integer such that $\gH^k\!\big(X,\CO(D(e))\big)=0$ for $e\in E^\times$ and all $k\geq k_0$, then also $\gH^k(X,\CE)=0$ for all $k\geq k_0$. In particular, this implies that if for all $e\in E^\times$ the line bundle
\begin{enumerate}
\item[{\rm (i)}] $\CO_X(D(e))$ is acyclic, then so is $\CE$. 

\smallskip

\item[{\rm (ii)}] $\CO_X(D(e))$ is immaculate, then so is $\CE$. 
\end{enumerate}
\end{theorem}

\begin{remark}
As the set of the divisors $D(e)$ is finite, Theorem~\ref{thm:VanThm2} is really Theorem~\ref{thm:VanThm} from the introduction.
\end{remark}

Remark~\ref{rem:LBPol} or Demazure vanishing immediately implies the

\begin{corollary}
\label{coro:Neflic}
A nefly decorated toric sheaf $\CE$ is acyclic. 
\end{corollary}

\begin{remark}
As a quotient of a direct sum of nef line bundles, a nefly decorated $\CE$ is actually also nef, cf.~\cite[Theorem 6.2.12]{lazarsfeld}.
\end{remark}

\section{Acyclicity of toric sheaves}
\subsection{Proof of Theorem~\ref{thm:Neflic}}
\label{subsect:Proof}
First, we may assume that without loss of generality $D=0$ for the twisting sheaf $\CO_X(D)$. If $D=\sum a_\rho D_\rho\not=0$, we pass to the toric sheaf $\CE':=\CE(D)$. This merely translates the Weil decoration by $D$, namely $\CD'(e)=\CD(e)+D$ for all $e\in E^\times$. In particular, $\CD'(\eta)=\CD(\eta)+D$ and $\widehat D'=\widehat D+D$ and thus $\mu_{\rho}'=a_{\rho}+\mu_{\rho}$ and $\lambda_{\rho}'=a_{\rho}+\lambda_{\rho}$, cf.\ \ssect{subsec:WeilDeco}. The assumption of Theorem~\ref{thm:Neflic} therefore becomes
\begin{equation}
\label{eq:FundEq}
\sum_{\rho\in\pc}\mu_{\rho}\geq\sum_{\rho\in \foc(\pc)(1)} \ffoc_\rho\cdot
\lambda_{\rho}
\end{equation}
for all primitive relations $\pc$ whose primitive inequality $\pRel(\pc)$
defines an extremal ray of $\NE(X)\subseteq\Cl(X)^\vee$. 

\medskip

We claim that this forces $\CE$ to be nefly decorated; actually, even $\widehat D$ becomes nef. Indeed, let
$\CD(e)=\sum b_\rho D_\rho$ be a divisor of the Weil decoration. Then
$\CD(e)\geq\CD(\eta)$ implies $\lambda_\rho\geq b_\rho\geq\mu_\rho$ and thus
\[
\langle\CD(e),\pRel(\pc)\rangle\;=\;\sum_{\rho\in\pc}b_{\rho}-\hspace{-0.7em}\sum_{\rho\in\foc(\pc)(1)}\hspace{-0.7em}\ffoc_\rho\cdot
b_{\rho}\;\geq\;\sum_{\rho\in\pc}\mu_{\rho}-\hspace{-0.7em}\sum_{\rho\in\foc(\pc)(1)}\hspace{-0.7em}\ffoc_\rho\cdot \lambda_{\rho}\;\geq\;0
\]
where $\langle\cdot\,,\cdot\rangle$ denotes the pairing in $\Z^{\Sigma(1)}$. Since this holds for all extremal primitive collections and also for $\widehat D$, we obtain that $\CD(e)$ and $\widehat D$ are nef.\hfill$\Box$

\subsection{Geometric interpretation of the assumption~\eqref{eq:PS}}
\label{subsect:GeoInt}
Though the Perlman-Smith assumption is too strong for acyclity, it is nevertheless worthwile giving a polyhedral interpretation. Towards that end, we continue to assume $D=0$ and write the resulting Equation~\eqref{eq:FundEq} as
\begin{equation}
\label{eq:FundEqB}
\fbox{$\langle\CD(\eta),\pRel(\pc)\rangle\;=\;\sum_{\rho\in\pc}\mu_{\rho}-\sum_{\rho\in\foc(\pc)(1)}\ffoc_\rho\cdot\mu_{\rho}$}\;\geq\;\fbox{$\sum_{\rho\in \foc(\pc)(1)}\ffoc_\rho\cdot
(\lambda_{\rho}-\mu_{\rho})$}. 
\end{equation}
Since the assumption of Theorem~\ref{thm:Neflic} implies that $\CE$ is nefly decorated, we may straight away assume that $\nabla(e)\in\Pol^+(\Sigma)$ for all $e\in E^\times$.

\subsubsection{The left hand side}
Since $\pRel(\pc)\in\Cl(X)^\vee$, the pairing only depends on the class $[\CD(\eta)]\in\Cl(X)$ and not on the concrete representative $\CD(\eta)$. 

\medskip

First, for every primitive collection giving rise to an extremal ray of the Mori cone, the primitive relation $\pRel(\pc)\in\NE(X)\subseteq\Cl(X)^\vee$ may be regarded as the class of a curve $[\ko{\orb}(\tau)]$ given by a wall $\tau\in \Sigma(d-1)$ of the fan. Hence
\[
\langle\CD(e),\pRel(\pc)\rangle=\langle[\CD(e)],\,\ko{\orb}(\tau)\rangle
\]
where the right hand side denotes the intersection product. 

\medskip

Second, each wall $\tau$ corresponds to an edge $k_\tau(\Delta)$ of any ample polytope $\Delta$ and more generally of a nef polytope if we allow edges to degenerate to a vertex. In particular, we obtain for every $e\in E^\times$ a (possibly degenerate) edge $k_\tau(e)\leq\nabla(e)$. 
Then $\langle[\CD(e)],\,\ko{\orb}(\tau)\rangle$ is precisely the
lattice length $\ell(k_\tau(e))$. In particular, the pairing of the left hand side is given by 
\[
\langle\CD(\eta),\pRel(\pc)\rangle=\ell(k_\tau(\eta))
\]
for $\pc=\pc(\tau)$.

\subsubsection{The right hand side}
By definition, the focus $\ffoc(\pc)=\sum_{\rho\in \foc(\pc)(1)} \ffoc_\rho\cdot e_\rho$ introduced in Section~\ref{sec:Intro} has the same image under $\pi\colon\Z^{\Sigma(1)}\to N$ as the incidence vector 
\[
e_{\pc}:=\sum_{\rho\in\pc}e_\rho 
\]
of the primitive collection $\pc\subseteq\Sigma(1)$. Next, we may regard $\widehat D-\CD(\eta)=\widehat\nabla-\nabla(\eta)$ as a formal difference of polytopes $\widehat\nabla\supseteq\nabla(\eta)$ in $\Pol(\Sigma)$. From~\eqref{eq:NablaFace} we obtain for $\sigma=\sigma(\pc)$ a difference vector $v(\sigma)$ in $M$ determined up to $\sigma^\perp$ pointing from $\face(\nabla(\eta),\sigma)$ to $\face(\widehat\nabla,\sigma)$, which we somehow abusively write
\[
v(\sigma):=\face(\widehat\nabla,\sigma)-\face(\nabla(\eta),\sigma)\in M/\sigma^\perp.
\]
Via the standard embedding of $M$ into $\Z^{\Sigma(1)}$, pairing $v(\sigma)$ with $\ffoc(\pc)$ yields
\begin{align*}
\sum_{\rho\in\foc(\pc)(1)}\big(\widehat D-\CD(\eta)\big)_\rho\cdot\ffoc_\rho&\;=\;-\langle v(\sigma),\; e_{\pc}\rangle\\[-10pt]
&\;=\;\min\langle\nabla(\eta),\,e_{\pc}\rangle-\min\langle\widehat\nabla,\,e_{\pc}\rangle
\end{align*}
and so magically eliminates the distinguished cone $\sigma$ of the primitive collection $\pc$. 

\medskip

Putting everything together, Assumption~\eqref{eq:PS} reads for $D=0$ and a nefly decorated $\CE$ in polyhedral terms as follows: For any extremal primitive collection $\pc=\pc(\tau)\subseteq\Sigma(1)$ associated with the wall $\tau\in\Sigma(d-1)$, the inequality
\begin{equation}
\label{eq:PSInter}
\ell\big(k_{\tau}(\eta)\big)\;\geq\;\min\langle\nabla(\eta),\,e_{\pc}\rangle-\min\langle \widehat\nabla,\,e_{\pc}\rangle
\end{equation}
holds true. Finally, the difference on the right hand side can be expressed as
\[
\sum_{\rho\in\pc}\big(\!\min\langle\nabla(\eta),\,\rho\rangle-\min\langle\widehat\nabla,\,\rho\rangle\big)
\]
where each summand encodes the distance between the $\rho$-facets of 
$\nabla(\eta)$ and $\widehat\nabla$.

\begin{example}
\label{exam:NefNotPS}
Equation~\eqref{eq:PSInter} always holds on a projective plane and more generally on any projective space $\P^n$, for there is only one primitive collection $\pc=\{\rho_0,\ldots,\rho_n\}$ with $e_\pc=0$. But already on the first Hirzebruch surface $\CH$, assumption~\eqref{eq:FundEqB} is strictly stronger than nefly decorated.

\medskip

Indeed, for the fan given in Figure~\ref{fig:NefNotPS} below, the Mori cone has two extremal rays $R_1$ and $R_2$ corresponding to the primitive collections $\pc_1=\{\rho_2,\,\rho_4\}$ and $\pc_2=\{\rho_1,\,\rho_3\}$, respectively. Now again $e_{\pc_1}=0$, but $e_{\pc_2}=e_{\rho_2}$. Since $R_2$ is induced by the ray $\rho_2$, we obtain the possibly nontrivial inequality 
\[
\ell(k_{\rho_2}(\eta)\;\geq\;\min\langle\nabla(\eta),\,\rho_2\rangle-\min\langle \widehat\nabla,\,\rho_2\rangle.
\]
Next, let $E=\C v\oplus\C v'$ with stratification $0$, $\S=\C v\setminus\{0\}$, $\S'=\C v'\setminus\{0\}$ and generic stratum $\eta$, which we decorate by the polytopes in Figure~\ref{fig:NefNotPS}, that is,
\[
0\mapsto\widehat\nabla,\quad\S\mapsto\nabla,\quad\S'\mapsto\nabla',\quad\eta\mapsto\nabla\cap\nabla'.
\]
This defines the Weil decoration of $\CO_\CH(\nabla)\oplus\CO_\CH(\nabla)'$ on $\CH$. Clearly, it is nefly (even amply) decorated. However,
\[
\min\langle\nabla(\eta),\,\rho_2\rangle-\min\langle \widehat\nabla,\,\rho_2\rangle=0-(-3)=3
\]
while $\ell(k_{\rho_2}(\eta)=2$, which violates Inequality~\eqref{eq:FundEqB}.

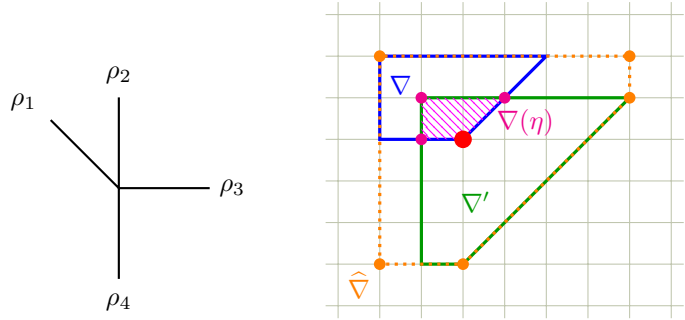
\begin{figure}[!ht]
\begin{tikzpicture}[scale=0.3]
\draw[thick, color=black]
  (0,0) -- (4,0) (0,-4) -- (0,4) (0,0) -- (-3.0,3.0);
\draw[thick, color=black]
  (5,0) node{$\rho_3$} (0,5) node{$\rho_2$} (-4.2,3.8) node{$\rho_1$} (0,-5) node{$\rho_4$};
\end{tikzpicture}
\hspace{20pt}
\begin{tikzpicture}[scale=0.55]
\draw[color=oliwkowy!40] (-1.3,-1.3) grid (7.3,6.3);
\draw[very thick, color=green]
  (1,0) -- (2,0) -- (6,4) -- (1,4) -- cycle;
\draw[thick, color=green]
  (2.3,1.5) node{$\nabla'$};
\draw[very thick, color=blue]
  (0,3) -- (2,3) -- (4,5) -- (0,5) -- cycle;
\draw[very thick, color=blue]
  (0.5,4.3) node{$\nabla$};
\draw[very thick, color=orange, dotted]
  (0,0) -- (0,5) -- (6,5) -- (6,4) -- (2,0) -- cycle; 
\fill[color=orange]
  (0,0) circle (4pt) (0,5) circle (4pt) (6,5) circle (4pt) (6,4) circle (4pt) (2,0) circle (4pt);
\draw[very thick, color=orange]
  (-0.5,-0.5) node{$\widehat\nabla$};
\fill[pattern color=magenta, pattern=north west lines]
  (2,3) -- (1,3) -- (1,4) -- (3,4) -- cycle;
\fill[color=magenta]
  (1,3) circle (4pt) (1,4) circle (4pt) (3,4) circle (4pt);
\draw[very thick, color=magenta]
  (3.5,3.4) node{$\nabla(\eta)$};
\fill[color=red]
  (2,3) circle (6pt);
\end{tikzpicture}
\caption{Left hand side: The fan of the first Hirzebruch surface. Right-hand side: The Weil decoration of $\CO_\CH(\nabla)\oplus\CO_\CH(\nabla)'$.}
\label{fig:NefNotPS}
\end{figure}
\end{example}

\subsection*{Declarations}

Ethical approval \& Data Availability: not applicable. \\
Funding: The second named author was partially supported by the project PRIN 2022K48YYP \emph{Unirationality, Hilbert schemes, and singularities}.\\

\end{document}